\documentclass[preprint,12pt]{elsarticle}

\usepackage[top=1in,bottom=1in,left=1in,right=1in]{geometry}

\usepackage{amsmath}
\usepackage{amssymb}
\usepackage{mathtools}
\usepackage{bm}

\usepackage{graphicx}
\usepackage{enumerate}
\usepackage{booktabs}
\usepackage{float}
\usepackage{algorithm}
\usepackage{algorithmic}

\usepackage{xcolor}
\usepackage{hyperref}
\hypersetup{colorlinks=true,linkcolor=blue,citecolor=blue,urlcolor=blue}

\newtheorem{thm}{Theorem}
\newtheorem{cor}[thm]{Corollary}
\newtheorem{lem}[thm]{Lemma}
\newtheorem{prop}[thm]{Proposition}

\newtheorem{exmp}[thm]{Example}
\newtheorem{rem}[thm]{Remark}

\allowdisplaybreaks[2]

\newcommand{\p}{\partial}

\newcommand{\pd}[2]{\frac{\partial #1}{\partial #2}}
\newcommand{\pdd}[2]{\frac{\partial^2 #1}{\partial #2^2}}
\newcommand{\R}{\mathbb{R}}
\newcommand{\E}{\mathbb{E}}
\newcommand{\G}{\mathcal{G}}

\newcommand{\T}[1]{\text{#1}}

\biboptions{sort&compress}

\myfooterfont{\normalfont\footnotesize}
\myfooter[C]{\thepage}

\begin{document}

\begin{frontmatter}

\title{Neural Operator Approximation for Yau--Yau Nonlinear Filtering: Stability, Approximation, and Fast Computation}

\author[CNU]{Xiaopei Jiao}\ead{jiaoxiaopei@cnu.edu.cn}
\author[RUC]{Xiuqiong Chen\corref{cor1}}\ead{cxq0828@ruc.edu.cn}
\author[THU]{Stephen S.-T. Yau\corref{cor1}}\ead{yau@uic.edu}

\cortext[cor1]{Corresponding authors.}

\address[CNU]{Capital Normal University, Beijing, China}
\address[RUC]{Renmin University of China, Beijing, China}
\address[THU]{Tsinghua University, Beijing, China}

\tnotetext[1]{This paper was not presented at any IFAC meeting.
This work was supported by grants of the National Natural Science Foundation
of China (No.~12501613, No.~12671549, No.~42450242) and the Tsinghua
University Education Foundation fund.}

\begin{abstract}
Since the introduction of the Kalman filter, nonlinear state estimation has become an important topic in areas such as aerospace engineering, sensor design, and target tracking. Within the optimal filtering framework described by the Duncan–Mortensen–Zakai (DMZ) equation, Yau–Yau filtering (YYF) employs an offline–online decomposition: the Fokker–Planck–Kolmogorov (FKE) equations, initialized by a prescribed set of basis functions, are solved offline, and the prior density is projected onto these basis functions during the online stage. However, this approach still faces two computational challenges: solving multiple FKEs offline can be computationally expensive, and the online projection becomes costly when filtering is performed repeatedly or in high-dimensional settings. To address these limitations, we develop a neural operator filter (NOF) that directly learns the mapping from the prior density to the posterior density, thereby eliminating both the expensive offline FKE solves and the online projection. The underlying operator is approximated using a Fourier neural operator (FNO). Numerical experiments on several challenging filtering systems verify both the effectiveness and accuracy of our method, and we further provide a rigorous theoretical analysis of its stability and approximation error.
\end{abstract}

\begin{keyword}
State estimation
\sep Neural operator
\sep Nonlinear filtering
\sep Forward Kolmogorov equation
\sep Yau--Yau filter
\end{keyword}

\end{frontmatter}

\section{Introduction}
Nonlinear filtering aims to estimate the internal state of a dynamic system from a sequence of noisy observations, extending the foundational Kalman filter \cite{kalman1960new,kalman1961new} to manage complex, nonlinear relationships. Its importance is demonstrated by transformative applications across numerous fields. In autonomous navigation, it enables precise real-time tracking by fusing data from GPS, inertial sensors, and LiDAR~\cite{mohamed1999adaptive}. For meteorology and oceanography, it is indispensable for data assimilation, dramatically improving the accuracy of weather and climate forecasts~\cite{evensen2003ensemble}. Throughout these diverse applications, nonlinear filtering serves as a critical enabling technology, driving advancement and providing deeper insights across science and engineering.

Current nonlinear filtering algorithms can be broadly organized into two categories. The first category involves finite-dimensional approximations of the posterior distribution, where algorithms approximate the posterior probability using finite-dimensional statistics, typically employing a Gaussian distribution or its weighted average as a proposal distribution, essentially relying on second-order moment approximations of the conditional density. Fundamental methods like the Extended Kalman Filter (EKF) and the Unscented Kalman Filter (UKF) and the Ensemble Kalman Filter (EnKF)~\cite{yau2024principles} fall into this group; however, they can yield significant approximation errors or even diverge for systems with strong nonlinearities or non-Gaussian, multi-modal densities. The second category comprises methods that directly fit the posterior distribution, primarily using Monte Carlo sampling, including Particle Filter (PF) and its feedback form~\cite{arulampalam2002tutorial,yang2013feedback}. While PF is widely applicable for non-Gaussian and highly nonlinear systems, it suffers from particle degeneracy and the curse of dimensionality. Recently, various deep learning architectures have been integrated into nonlinear filtering~\cite{qian2022deep,chen2024uniform,shi2024dglg,dong2026difnet}.

In the optimal filtering, the least squares estimate for a nonlinear system's state is the conditional expectation given the observation history, meaning the (conditional) density function contains all statistical information for optimal estimation. In the 1960s, Duncan-Mortensen-Zakai (DMZ) equation~\cite{duncan1967probability,mortensen1966optimal,zakai1969on} as a stochastic differential equation was proposed to describe the conditional density. It faces challenges to be solved accurately. Based on this, Yau and Yau proposed the frozen coefficient method, introduced a class of parabolic PDE, and took a rigorous convergence analysis, named the Yau-Yau filter~\cite{yau2000real,yau2008real}, which decomposes the state estimation into online and offline parts.  This two-phase procedure offered several advantages, such as the real-time merit, memoryless structure, and strict convergence guarantees. Following on this, a number of generalized algorithms and extensions were developed \cite{Luo_Complete2013,li2022solving,jiao2026extended}, etc.

Real-time response in nonlinear filtering is an important criterion in many practical scenarios, such as navigation, autonomous vehicles, and missile guidance. In recent years, there has been a breakthrough in learning nonlinear functional mappings, called neural operators such as Fourier neural operator (FNO) \cite{li2021fourier}, deep operator network (DeepONet) \cite{lu2021learning}, graph kernel network, etc. In this paper, we shall consider the general nonlinear time-invariant filtering systems and design a novel neural operator filter (NOF) by integrating Fourier neural operator (FNO) and YYF to accelerate posterior prediction of the Forward Kolmogorov equation (FKE). Our motivation is to learn the mapping between the densities in infinite dimensional functional space so that the posterior density can be immediately obtained when receiving the upcoming new observation. The datasets consists of two parts including Monte Carlo simulations of the true dynamics and data generated from a Gaussian random field. 

The contributions of our paper are summarized point by point as follows.
\begin{itemize}
\item For the FKE equation, we proved the stability of its associated forward operator (Theorem \ref{Thm_2}) and verified the existence of a universal approximating neural network (Corollary \ref{Cor:FKE_verification}).
\item Following the established FNO approximation (Lemma \ref{sec4.3:lemma1}), we integrate single-step forward operator error into iterative YYF density estimation (Theorem \ref{sec4.3_thm1}). 
\item By utilizing FNO, NOF significantly speed up the repetitive computation in YYF off-line stage before observation arrives.
\item For low- and medium-dimensional nonlinear systems, our NOF provides competitive tracking accuracy with reference methods.
\end{itemize}

The paper will be organized as follows. Section 2 will introduce basic filtering problem and formulation. Section 3 presents the methodology of our new method, the neural network, training process and the algorithm. In Section 4, stability analysis, the approximation error bounds will be given. Section 5 will show numerical results of benchmark filters. We conclude the paper in Section 6 and postpone the proof details to the Appendix.

\textbf{Mathematical Notations:}
Let $\R^n$ be Euclidean $n$-space, $\Omega\subset\R^n$ bounded with boundary $\partial\Omega$. For a domain $U$, $C^\infty(U)$, $C^r(U)$, $L^2(U)$, $L^\infty(U)$, $H^r(U)$, $W^{\ell,q}(U)$ denote the spaces of smooth, $r$-times continuously differentiable, square-integrable, essentially bounded functions, and Sobolev spaces with weak derivatives up to order $r$ in $L^2$ and up to $\ell$ in $L^q$, respectively. For smooth $\phi$, set $\nabla\phi=(\partial_1\phi,\dots,\partial_n\phi)^\top$, $\Delta\phi=\sum_i\partial_i^2\phi$ with $\partial_i:=\partial/\partial x_i$, and let $\nabla\cdot$ be divergence; $\mathbb{E}$ denotes expectation, $I_n$ the $n\times n$ identity.

\section{Basic concepts and Preliminaries}\label{sec:preliminary}
%

In this paper, we shall consider the \textbf{nonlinear time-invariant filtering} (NTIF) system:
\begin{equation}\label{NFsystem}
    \left\{
    \begin{alignedat}{2}
        d x(t) =& f(x(t))dt + \Gamma dw(t),\quad &&x(0)\sim p_0(x)\\
    d y(t) =& h(x(t))dt + dv(t), \quad&&y(0) = 0
    \end{alignedat}
    \right.
\end{equation}
where $x(t)\in\R^n$ and $y(t)\in\R^m$ denote system state and observation process. The function $f:\R^n\to\R^n$ is the drift and $h:\R^n\to\R^m$ is the observation function. $\{w(t)\}_{t\ge0}$ and $\{v(t)\}_{t\ge0}$ represent independent Brownian motion with covariance matrices $\E[dw(t)dw(t)^\top] = Q$ and $\E[dv(t)dv(t)^\top] = S$ which are both independent of initial state $x(0)$. Let $\Gamma Q \Gamma^\top = I_n$ be an identity matrix.

Given all observations up to time $t$, i.e., $\mathcal{Y}_t:=\{y_s:0\le s\le t\}$, the conditional density function (CDF) carries all the information needed for optimal estimation in the mean-square sense. Its unnormalized version satisfies the Duncan--Mortensen--Zakai (DMZ) equation:
\begin{equation}\label{eq-dmz}
    \left \{
    \begin{aligned}
    d\sigma(t, x) & =  \mathcal{L}\sigma(t, x)dt + \sigma(t, x)h(x)^\top S^{-1}\circ dy(t), \\
    \sigma(x,0)  & =  \sigma_0(x),
    \end{aligned}
    \right.
\end{equation}
where the operator $\mathcal{L}(\cdot)$ is defined as follows:
\begin{equation}\label{operator-L}
    \mathcal{L}(\phi) := \frac{1}{2}\sum_{i=1}^n \frac{\partial^2 \phi}{\partial x_i^2} - 
    \sum_{i=1}^n\frac{\p f_i}{\p x_i} \phi - \sum_{i=1}^nf_i\frac{\p\phi}{\p x_i}
\end{equation}
and $\sigma_0(x)$ is the unnormalized density of the initial state $x_0$. Once the solution is obtained, the CDF and the state estimate (the conditional expectation) follow as
\begin{equation}
p(t, x) = \frac{\sigma(t, x)}{\int_{\R^n} \sigma(t,x)dx}\quad\text{and}\quad\hat{x}(t) = \frac{\int_{\R^n} x\sigma(t, x)dx}{\int_{\R^n} \sigma(t,x)dx}
\end{equation}
In order to get a practical algorithm in the sense ``path-wise'', a robust DMZ equation was proposed by taking the corresponding Davis transformation \cite{davis1980multiplicative}: 
\[
u(x,t) = \exp(-h(x)^\top S^{-1}y_t)\,\sigma(x,t)
\]
and
\begin{equation}\label{rdmz}
    \left \{
    \begin{aligned}
    \frac{\partial u}{\partial t}(x,t) =&  \frac{1}{2}\Delta u(t,x) + (-f(x) + \nabla K(t,x))\cdot\nabla u(t,x)\\
    &+\Big(-\nabla\cdot f(x) - \frac{1}{2}|h(x)|^2 + \frac{1}{2}\Delta K(t,x)\\
    &-f(x)\cdot \nabla K(t,x) + \frac{1}{2}|\nabla K(t,x)|^2\Big)u(t,x)\\
    u(x,0) =& \sigma_0(x).  \\
    \end{aligned}
    \right.
\end{equation}
where $K(t,x)=y(t)\cdot h(x)$. On a time partition $P_N=\{0=\tau_0<\tau_1<\cdots<\tau_N=T\}$ with $\Delta\tau=\tau_i-\tau_{i-1}$, although the robust DMZ equation handles a single observation path, it can be further separated into ``off-line'' and ``on-line'' parts, which yields the so-called Forward Kolmogorov equation (FKE). For $t\in[\tau_{i-1},\tau_i]$, 
\begin{equation}\label{FKE_filter}
\left\{
\begin{aligned}
\frac{\partial \tilde{u}_i}{\partial t}(t,x)
&= \frac12 \Delta \tilde{u}_i(t,x)
   + F(x)\cdot \nabla \tilde{u}_i(t,x) + J(x)\tilde{u}_i(t,x),\quad 
  \\[6pt]
\tilde{u}_i(\tau_{i-1},x)
&= \exp\!\left(
      \big( y(\tau_{i-1}) - y(\tau_{i-2}) \big)\cdot h(x)
    \right)
    \tilde{u}_{i-1}(\tau_{i-1},x).
\end{aligned}
\right.
\end{equation}
where $F(x)=-f(x)$ and $J(x) = -(\nabla\cdot f(x) + \frac{1}{2}\sum_{j=1}^m h_j^2(x))$.  $\tilde{u}_{i}(\tau_i, x)$ is obtained by
\[
 \tilde{u}_i(\tau_i, x)
= \exp\!\left(
     \sum_{j=1}^m y_j(\tau_{i-1})\, h_j(x)
  \right)u_i(\tau_i, x)
 .
\]
where $u_i(t,x)$ is a piece-wise truncated solution of $u$ on $[\tau_{i-1},\tau_i]$ obtained by setting $y(t) = y(\tau_{i-1})$.
\begin{rem}
To simplify the DMZ equation, we assume that $\Gamma Q\Gamma^\top=I_n$, which gives the diffusion term $\frac12\Delta\tilde u_i(t,x)$ in the FKE. This does not restrict our algorithm, since the FKE contains no explicit time dependence under our setting. Throughout the theoretical analysis we also set $S=I$ without loss of generality; the numerical examples in Section 5 illustrate the general case $S\ne I$.
\end{rem}


\section{Neural operator-based filter algorithm}
\subsection{Neural operator}\label{sec::3.1}
We shall consider a parametric partial differential equation formulated as follows.
\begin{equation}
\left\{
\begin{alignedat}{2}
&\mathcal{L}_a[u](x) = 0, \quad&&\text{in}\quad x\in\Omega\\
&u(x) = g(x)\quad&&\text{in}\quad x\in\p\Omega
\end{alignedat}
\right.
\end{equation}
where $\mathcal{L}$ is a general (nonlinear) differential operator parameterized by $a\in\mathcal{A}$, $\Omega\subset\R^n$ is the computational domain with boundary $\partial\Omega$, $u\in\mathcal{U}$ is the solution in a Banach space, and $g$ is the boundary condition. A typical problem is then to approximate the solution operator $\G:\mathcal{A}\to\mathcal{U}$, $\G(a)=u$.

Given a set of input--output pairs $\{a_i,u_i\}_{i=1}^N$ with $a_j\sim\mu$, the goal is to find a neural network $\G_{\theta^*} : \mathcal{A}\to\mathcal{U}$ such that $\G_{\theta^*}\approx\G$, i.e., to minimize a suitable operator distance $\min_\theta \mathbb{E}_{a\sim\mu}[\mathrm{dist}(\G_\theta(\cdot),\G(\cdot))]$.

Inspired by the Green-function (convolutional) representation of the solution, Li et al.~\cite{li2021fourier} designed $\G_\theta$ as 
\[
\G_\theta := \mathcal{Q}\circ \mathcal{L}_L\circ \mathcal{L}_{L-1}\circ \cdots \mathcal{L}_{1}\circ \mathcal{P}
\]
where $\mathcal{P},\mathcal{Q}$ are the encoder and decoder, and the kernel integral operator $\mathcal{L}_{i}$ is
\[
(\mathcal{L}_{i}v)(x) = \int_{\Omega}\kappa_{\phi}(x, y) v(y)dy
\]
with a kernel $\kappa_\phi$ learned from data. To filter high-frequency components and accelerate the computation, the Fourier integral operator is used instead:
\[
(\mathcal{L}_{i}v)(x) = \mathcal{F}^{-1}(R_\phi\cdot \mathcal{F}v)(x)
\]
where $\mathcal{F}$ and $\mathcal{F}^{-1}$ denote the Fourier transform and its inverse, $R_\phi$ is a learned kernel, and both are truncated to the maximal Fourier modes $k_{\T{max}}=\{k\in\mathbb{Z}^n: |k_j|\le k_{\text{max},j},\, j=1,\dots,n\}$.

\subsection{Methodology for neural operator filter}\label{sec::3.2}
The FKE plays a central role in propagating the CDF. Restricting the solution to the ball $B_R\subset\R^n$, Eq.~(\ref{FKE_filter}) becomes the following initial-boundary value problem.
\begin{equation}\label{FKE_eq}
\left\{
\begin{aligned}
\frac{\partial \tilde{u}}{\partial t}(t,x)
=& \frac12 \Delta \tilde{u}(t,x)
   + F(x)\cdot \nabla \tilde{u}(t,x)\\
   &+ J(x)\tilde{u}(t,x),&&\text{on  } [0, \Delta \tau]\times B_R
  \\
\tilde{u}(0,x) =& \phi(x), &&\text{on  }\{t=0\}\times B_R\\
\tilde{u}(t, x) =& 0, &&\text{on  }[0,\Delta \tau]\times\p B_R
\end{aligned}
\right.
\end{equation}
where
\begin{equation}\label{dfn_F_J}
F(x) := -f(x),\quad J(x) := -\bigg[\nabla\cdot f(x) + \frac{1}{2}\sum_{j=1}^m h_j^2(x)\bigg]
\end{equation}
Given a filtering system (\ref{NFsystem}), the FKE in Eq.~(\ref{FKE_filter}) can be solved iteratively to get $\{\tilde{u}_i(t, x)\}_{i\ge 1}$ at interval $[\tau_{i-1}, \tau_i]\times B_R$ by considering $\tilde{u}_i(\tau_{i-1}, x)$ as initial. Therefore, we shall consider the forward operator of FKE, denoted by $\G^\dagger$ in $x\in B_R$.
\begin{equation}\label{FKE_forward_map}
\G^\dagger: L^2(B_R)\to L^2(B_R),\quad\tilde{u}_i(\tau_{i-1}, x)\mapsto\tilde{u}_i(\tau_{i}, x)
\end{equation}
for $i=1,\dots,N$, or equivalently $\phi(x)\mapsto\tilde u(\Delta\tau,x)$, since Eq.~(\ref{FKE_eq}) is autonomous. Because the system is time-invariant ($f,h$ do not depend on time), $\G^\dagger$ is independent of $i$.  

\subsubsection{Training strategy}
We design an operator network $\G_\theta$ to approximate the mapping from the prior CDF $\tilde u_i(\tau_{i-1},x)$ to the posterior CDF $\tilde u_i(\tau_i,x)$, where $\theta=\{W_i,b_i\}$ collects the trainable parameters. As in Section 3.1, the architecture is the composition $\G_\theta=\mathcal{Q}\circ\mathcal{L}_L\circ\cdots\circ\mathcal{L}_1\circ\mathcal{P}$, where $\mathcal{P}:\R^{d_a}\to\R^{d_1}$ and $\mathcal{Q}:\R^{d_L}\to\R^{d_u}$ are the encoder and decoder that lift and recover the channel dimension ($d_a,d_u$ being the codomain dimensions of the input and output functions), and the Fourier convolution layer $\mathcal{L}_\ell$:
\[
\mathcal{L}_\ell(v_\ell) = \sigma[W_l v_\ell + b_l + \mathcal{F}^{-1}_N(\mathcal{R}_l\cdot \mathcal{F}_N(v_\ell))]
\]
\textbf{Dataset preparation.}
We build a multi-source dataset. First, random filter trajectories are sampled by Monte Carlo simulation to obtain pairs of prior and posterior CDFs $\{\tilde u_i^{(j)}(\tau_{i-1},x_k),\, \tilde u_i^{(j)}(\tau_i,x_k)\}_{i,j,k\ge1}$ where the posterior is computed by solving the Yau--Yau filter with a finite-difference method, $j$ indexes Monte Carlo trajectories, and $\{x_k\}\subset B_R$ are $N_c$ collocation points. This part of the data captures the intrinsic conditional-density evolution, so that the operator network learns the sufficient statistics of the filtering system; the resulting tuples of training dataset of NOF are prepared as $\{(x_k,\tilde u_i^{(j)}(\tau_{i-1},x),\tilde u_i^{(j)}(\tau_i,x_k))\}$.

Second, to improve generalization, we add FKE input--output pairs in which the initial condition is generated by zero-mean Gaussian random field (GRF), $\phi_j(x)\sim\mathcal{GP}(0,k_l(x,y))$, with covariance kernel $k_l(x,y)=\exp(-\|x-y\|^2/2l^2)$ and length scale $l>0$; the reference solutions are computed by a standard finite-difference method. This yields a second dataset of triplets $\{x_k,\phi_j(x),\tilde u_j(\Delta\tau,x_k)\}$, with $N_r$ sampled inputs and $N_c$ collocation points.

\textbf{Loss function.}
Our constructed operator network has the following loss function,
\begin{equation}
\begin{aligned}
&\mathcal{L}(\theta) = \mathcal{L}_{\T{MC}}(\theta) + \mathcal{L}_{\T{GP}}(\theta)\\
&= \frac{1}{N N_m N_c}\sum_{i=1}^{N}\sum_{j=1}^{N_m}\sum_{k=1}^{N_c}\underbrace{\bigg|\tilde{u}_i^{(j)}(\tau_{i}, x_k) - \G_\theta(\tilde{u}_i^{(j)}(\tau_{i-1}, \cdot))(x_k)\bigg|^2}_{\T{Monte-Carlo filtering}}\\
&+ \frac{1}{N_r N_c}\sum_{j=1}^{N_r}\sum_{k=1}^{N_c}\underbrace{\bigg|\tilde u_j(\Delta\tau, x_k) - \G_\theta(\phi_j(x))(x_k)\bigg|^2}_{\T{Gaussian random fields}}
\end{aligned}
\end{equation}
where $\theta$ collects the weights, biases, and truncated Fourier kernel parameters, and $\mathcal{L}_{\T{MC}}$, $\mathcal{L}_{\T{GP}}$ correspond to the Monte Carlo filtering dataset and the Gaussian random field dataset, respectively.

\subsubsection{Deployment in on-line filtering stage}
\begin{figure}[!htb]
    \centering
    \includegraphics[width=0.8\linewidth]{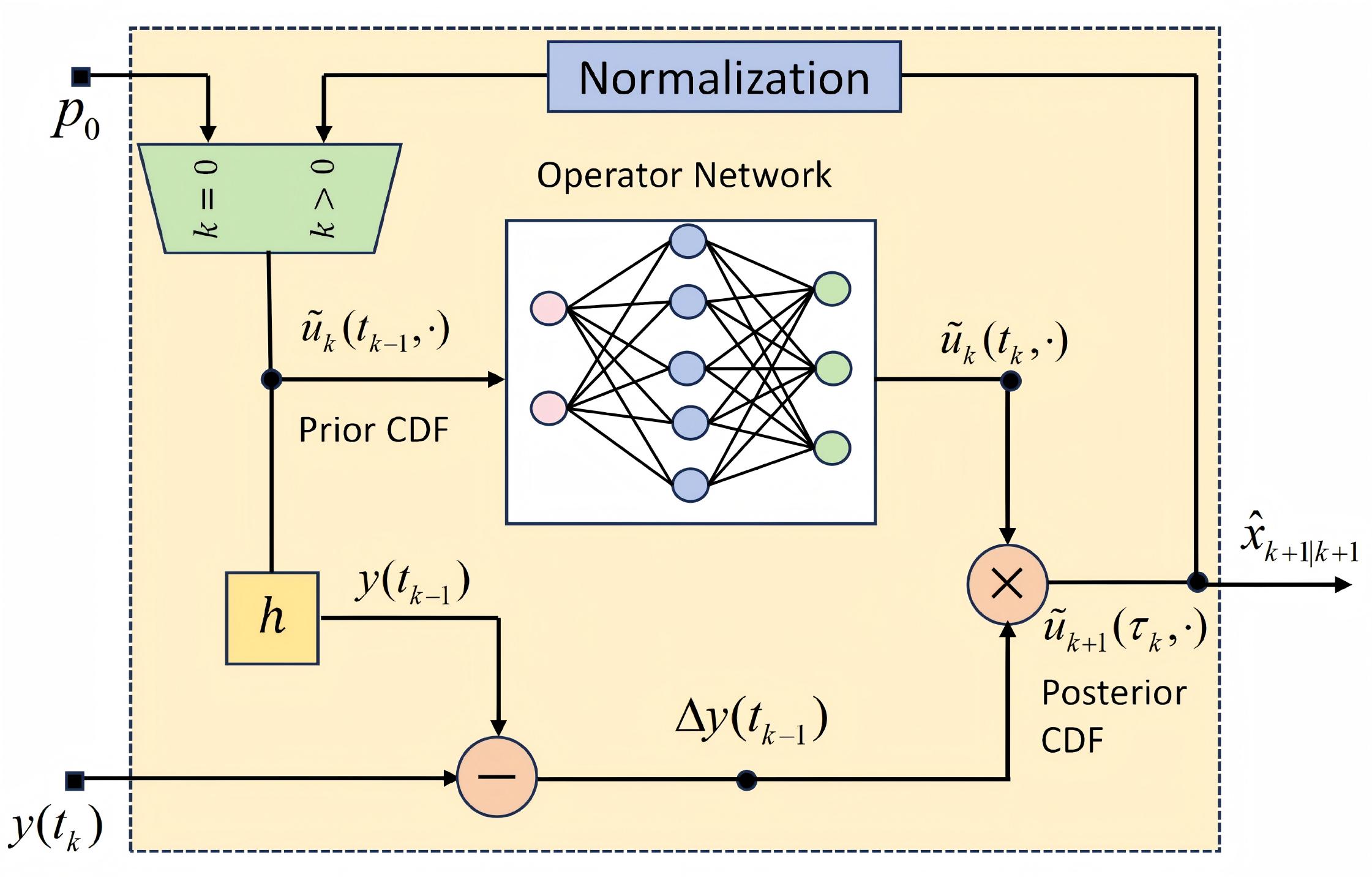}
    \caption{Block diagram of the neural operator filter.}
    \label{fig:pipeline}
\end{figure}

After training, we obtain a surrogate $\G_{\theta^*}\approx\G^\dagger$ of the forward operator, which can be deployed as new observations of the same system arrive. The online filtering process then combines this trained operator with the Yau--Yau filtering algorithm, as summarized in Algorithm \ref{NOF_pseudo_code} and illustrated in Fig.~\ref{fig:pipeline}.

%

\begin{algorithm}\small
\caption{Deployed procedure for NOF}\label{NOF_pseudo_code}
\label{alg:nof}
\textbf{Input:} Filter system $(f, h, \Gamma, Q, S)$ and spatial ball $B_R$; Time interval partition $\{0=\tau_0<\tau_1<\cdots<\tau_{N}=T\}$ with equal partition $\tau_i=\frac{T}{N}$; Initial distribution $p_0$; Trained operator network $\G_\theta$;\\
\textbf{Output:} Normalized approximated CDF $p(\tau_i, x)$; State estimation $\hat{x}(\tau_{k})$.
\begin{algorithmic}[1]
\STATE Initialize density $\tilde{u}_1(0,\cdot)=p_0(x)$
\FOR{$ k = 1 \rightarrow N $}
    \STATE Discretize on spatial grid $\{x_j\}$ and normalize prior CDF $\tilde{u}_k(\tau_{k-1}, x_j)\leftarrow \tilde{u}_k(\tau_{k-1}, x_j)/\sum_{j=1}^{N_j}\tilde{u}_k(\tau_{k-1}, x_j)$
    \STATE Compute posterior density of prediction by trained network $\tilde{u}_k(\tau_{k}, x_j)\leftarrow \G_\theta(\tilde{u}_k(\tau_{k-1}, \cdot))(x_j)$.
    \STATE  Update density by using latest observation and store posterior CDF $\tilde{u}_{k+1}(\tau_{k}, x_j) \leftarrow \exp(h(x_j)^\top S^{-1}(y_{\tau_k} - y_{\tau_{k-1}}))\cdot \tilde{u}_{k}(\tau_{k}, x_j)$.
    \STATE Compute the state estimate from the posterior density;\\
 $\hat{x}(\tau_{k})\leftarrow \sum_{j=1}^{N_j}\tilde{u}_{k+1}(\tau_{k}, x_j)x_j/\sum_{j=1}^{N_j}\tilde{u}_{k+1}(\tau_{k}, x_j)$.
\ENDFOR
\end{algorithmic}
\end{algorithm}

Compared with existing deep-learning filters, NOF is distinct: sequence and deep state-space models \cite{krishnan2015deep,chen2024uniform,qian2022deep} return only a conditional-mean estimate; KalmanNet \cite{revach2022kalmannet} inherits the Gaussian/linear assumptions of the Kalman recursion; and projection-based filters such as DGLG \cite{shi2024dglg} pay a basis-selection and re-projection cost at every observation. In contrast, NOF learns, once and offline, a surrogate $\mathcal{G}_\theta$ of the FKE forward map $\mathcal{G}^\dagger$, and reuses it online inside the Yau--Yau loop, propagating the full unnormalized conditional density without basis functions or online projection---the fixed-$N$ regime of Theorem \ref{sec4.3_thm1}.

\section{Theory}
In this section we analyze the approximation error of the NOF for the solution of robust DMZ equation, namely, the unnormalized conditional density. First, we establish the stability of the FKE forward operator under perturbations of the initial condition and the forcing term via an energy estimate. Second, we study the neural-network approximation of the FKE, both for the whole-space solution and for its truncation to a bounded ball $B_R$; the truncation error is made arbitrarily small by enlarging the ball radius. Finally, we analyze the approximation ability of the NOF on filtering process by providing an effective upper bound for a fixed time partition resolution.

\subsection{Stability of the forward operator Kolmogorov equation}
We first discuss the stability of the FKE (Eq.~(\ref{FKE_eq})): the solution should be insensitive, to some extent, to perturbations. We model perturbations of the forcing term and the initial condition by $R_i[\cdot]$ and $\delta\phi(x)$, giving the perturbed FKE
\begin{equation}\label{FKE_perturbation}
\left\{
\begin{alignedat}{2}
\frac{\partial \tilde{v}}{\partial t}(t,x)
=& \frac12 \Delta \tilde{v}(t,x)
   + F(x)\cdot \nabla \tilde{v}(t,x)\\ 
   &+ J(x)\tilde{v}(t,x) + R_i[\tilde{v}](t,x),\quad&&\text{  on  } [0, \Delta \tau]\times B_R
  \\
\tilde{v}(0,x)=& \phi(x)+\delta\phi(x),\quad &&\text{  on  }\{t=0\}\times B_R\\[0.5ex]
\tilde{v}(t, x) =& 0, \quad&&\text{  on  }[0,\Delta \tau]\times\p B_R
\end{alignedat}
\right.
\end{equation}
This analysis quantifies how the terminal error $|\tilde u-\tilde v|$ propagates starting from the initial error $|\delta\phi|$. Subtracting (\ref{FKE_eq}) from (\ref{FKE_perturbation}), the error $w(t,x)=\tilde v(t,x)-\tilde u(t,x)$ satisfies
\begin{equation}\label{Eq::w_R}
\left\{
\begin{alignedat}{2}
\pd{w}{t}(t,x) =& \frac{1}{2}\Delta w + F(x)\cdot \nabla{w} \\
&+ J(x)w + R_i[\tilde{v}]\quad&&\text{  on  } [0, \Delta \tau]\times B_R\\
w(0, x) =& \delta\phi(x),\quad &&\text{  on  }\{t=0\}\times B_R\\[0.5ex]
w(t, x) =& 0, &&\text{  on  }[0,\Delta\tau]\times\p B_R
\end{alignedat}
\right.
\end{equation}
The following proposition shows that the terminal error $w(T,\cdot)$ and its gradient are bounded by the initial perturbation $\delta\phi$ and the PDE residual $R_i[\tilde v]$; the proof is given in Appendix \ref{Proof_Prop4.1}.
\begin{prop}\label{Thm_1}
Let $w$ be a classical error solution satisfying Eq.(\ref{Eq::w_R}) and assume that $F(x)\in C^1([0, T]\times B_R)$ and $J(x)\in C([0, T]\times B_R)$. The error estimate will be given as
\begin{equation}\label{sec4.1::prop1}
\begin{aligned}
&\int_{B_R}|w(T,x)|^2 + \int_0^T\int_{B_R}|\nabla w(s,x)|^2 \le \int_{B_R}|\delta\phi(x)|^2 \\
&+ \sup_{x\in B_R}\bigg[-\nabla\cdot F(x) + 2J(x) + 1\bigg]\int_0^T\int_{B_R}|w(s,x)|^2 \\
&+ \int_0^T\int_{B_R}R_i[\tilde v]^2.
\end{aligned}
\end{equation}
\end{prop}
The preceding estimate shows that the forward map $\mathcal{G}^\dagger$ in (\ref{FKE_forward_map}) is Lipschitz, so no singularity of the posterior CDF $\tilde u(T,\cdot)$ can arise from a bounded prior $\tilde u(0,\cdot)$. Eq.~(\ref{sec4.1::prop1}) is an integrated form; its differential form is 
\begin{equation}\label{Diff_form_w}
\begin{aligned}
&\frac{1}{2}\frac{d}{dt}\int_{B_R}|w|^2 + \int_{B_R}\frac{1}{2}|\nabla w|^2\\
&\le \sup_{x\in B_R}\bigg[-\frac{1}{2}\nabla\cdot F(x) + J(x) + \frac{1}{2}\bigg]\int_{B_R}w^2 + \frac{1}{2}\int_{B_R}R_i[\tilde{v}]^2
\end{aligned}
\end{equation}
Setting $R_i[\tilde v]=0$ in (\ref{Diff_form_w}) yields the differential-form estimate for the FKE itself:
\begin{equation}\label{Diff_form_FKE}
\begin{aligned}
&\frac{1}{2}\frac{d}{dt}\int_{B_R}|\tilde u|^2 + \int_{B_R}\frac{1}{2}|\nabla \tilde u|^2\\ 
\le& \sup_{x\in B_R}\bigg[-\frac{1}{2}\nabla\cdot F(x) + J(x) + \frac{1}{2}\bigg]\int_{B_R}\tilde u^2 
\end{aligned}
\end{equation}
Based on Eq.(\ref{Diff_form_FKE}) and the Gronwall inequality technique \cite{evans2010partial}, we obtain

\begin{thm}\label{Thm_2}
Functional mapping $\mathcal{G}: \tilde{u}(0, x)\mapsto \tilde{u}(T, x)$ associated with the FKE has the Lipschitz property, i.e.,
\begin{equation}
\|\tilde{u}(T,\cdot)\|_{L^2(B_R)}\le L_\mathcal{G}\|\tilde{u}(0,\cdot)\|_{L^2(B_R)}
\end{equation}
with Lipschitz constant 
\begin{equation}\label{Lip_const}
L_\mathcal{G} = \exp(T\sup_{x\in B_R}\left[-\frac{1}{2}\nabla\cdot F(x) + J(x) + \frac{1}{2}\right]).
\end{equation}
\end{thm}
Detailed proof can be found in section \ref{Proof_thm4.1}.


\subsection{Neural-network approximation of the Forward Kolmogorov equation}\label{sec::4.2}
We shall formulate the FKE on the whole space domain $\R^n$ as follows, i.e., the Cauchy initial value problem.
\begin{equation}\label{global_FKE}
\left\{
\begin{aligned}
\frac{\partial \tilde u_g}{\partial t}(t,x)
=& \frac12 \Delta \tilde u_g(t,x)
   + F(x)\cdot \nabla \tilde u_g(t,x)\\
   &+ J(x)\tilde u_g(t,x),&&\text{on  } [0, \Delta \tau]\times \R^n
  \\
\tilde u_g(0,x) =& \phi(x), &&\text{on  }\{t=0\}\times \R^n
\end{aligned}
\right.
\end{equation}

We seek a network that approximates the solution of this time-dependent PDE at a fixed time: for the global FKE, we construct $\mathcal{U}^\varepsilon:\mathcal{X}\times[0,T]\to\mathcal{Y}$ mapping $(\phi,t)$ to $\tilde u_g$ with accuracy $\varepsilon$ (proof in Appendix \ref{pf:thm:FKE_verification}).

\begin{thm} \label{thm:FKE_verification}
Let $F$ be globally Lipschitz with constant $L_F$;
$\nabla\cdot F \in L^\infty(\mathbb{R}^n)$; $J$ is locally Lipschitz and bounded above by $\bar{J}<\infty$. Then, for any $\varepsilon>0$, $t\in[0,T]$, $\ell\in\mathbb{N}$, $q\ge1$, $\phi\in C^\ell$, and sufficiently large $R>0$, there exists a neural network
$\mathcal{U}^\varepsilon(\phi, t) : B_R \to \mathbb{R}$ such that
\begin{equation}\label{sec4.2:thm1:eq}
\begin{aligned}
\| \mathcal{U}^\varepsilon(\phi, t) - \tilde{u}_g(t, \cdot) \|_{L^q(B_R)}
&\le \varepsilon,\\
\ \text{and}\ \max_{t \in [0,T]} \| \mathcal{U}^\varepsilon(\phi, t) \|_{W^{\ell, q}(B_R)} 
&\le C_{\varepsilon,\ell}^B
\end{aligned}
\end{equation}
where $C_{\varepsilon,\ell}^B$ depends on $\varepsilon$, $R$, $\|\phi\|_{C^\ell}$, $L_F$, $\bar{J}$, and $T$.
\end{thm}
We next construct a neural-network approximation for the bounded-domain solution $\tilde u$ by combining the network error (\ref{sec4.2:thm1:eq}) with the gap $|\tilde u_g-\tilde u|$ (proof in Appendix \ref{pf:Cor:FKE_verification}).
\begin{cor}\label{Cor:FKE_verification}
Under the assumptions of Theorem \ref{thm:FKE_verification}, for any $\varepsilon>0$, $t\in[0,T]$, sufficiently large $R>0$, and $\phi$ satisfying the decay condition $\left(\int_{\mathbb{R}^n}\mathrm{e}^{\sqrt{1+|x|^2}}|\phi(x)|\,dx\right)^{1/2}<\infty$, there exists a neural network
$\mathcal{U}^\varepsilon(\phi, t) : B_R \to \mathbb{R}$ such that
\begin{equation}
\| \mathcal{U}^\varepsilon(\phi, t) - \tilde u(t, \cdot) \|_{L^2(B_R)}
\le \varepsilon + C_{\mathrm{trunc}}(R),
\end{equation}
where
\[
C_{\mathrm{trunc}}(R) \le C_{n,T} \mathrm{e}^{-\frac{9}{16}R}
\left( \int_{\mathbb{R}^n} \mathrm{e}^{\sqrt{1+|x|^2}} |\phi(x)| \, dx \right)^{1/2}
\]
is the truncation error of whole space solution $\tilde u_g$ truncated on $B_R$.
\end{cor}

\subsection{Approximation error of neural operator filter}


In the preceding sections we established that the required neural-network approximator exists (Corollary \ref{Cor:FKE_verification}). We now lift this pointwise-in-time approximation to a full operator approximation. Since the analysis of Fourier-based architectures is usually performed on a periodic domain, we embed $B_R$ into a torus $\mathbb{T}^n_L=[-L/2,L/2]^n\supseteq B_R$ of larger size $L$ with periodic boundary conditions, and introduce the following operators:
\begin{itemize}
    \item Extension operator $\mathcal{E}:C^r(B_R)\to C^r(\mathbb{T}^n_L)$ with $\mathcal{E}\phi|_{B_R}=\phi$, $\|\mathcal{E}\phi\|_{C^r(\mathbb{T}^n_L)}\le C_{\T{ext}}^{(C)}\|\phi\|_{C^r(B_R)}$, and $\|\mathcal{E}\phi\|_{H^r(\mathbb{T}^n_L)}\le C_{\T{ext}}^{(H)}\|\phi\|_{H^r(B_R)}$;
    \item Restriction operator $\mathcal{R}:C^r(\mathbb{T}^n_L)\to C^r(B_R)$, $\mathcal{R}f:=f|_{B_R}$;
    \item Trigonometric projection $P_N:=\mathcal Q_N\circ \mathcal E_N$.
\end{itemize}

We define the auxiliary operator $\hat\G(\phi):=\mathcal{U}^\varepsilon\circ\mathcal{R}\circ P_N\circ\mathcal{E}(\phi)$ and construct a FNO $\G_\theta(\phi):=\mathcal{R}\circ P_N\circ\mathcal{E}\circ\hat\G$ for $\phi\in C^r(B_R)$, and analyze the error between $\G^\dagger$ and $\G_\theta$ (proof in Appendix \ref{pf:sec4.3:lemma1}).

\begin{lem} \label{sec4.3:lemma1}
Let $\mathcal{G}^\dagger$ be the forward operator of FKE defined in \eqref{FKE_forward_map} on a time interval of length $\Delta\tau$, acting on $C^r(B_R)$. Under the assumptions of Theorem \ref{thm:FKE_verification}, consider the functional space $S_\phi=\{\phi: \left(\int_{\mathbb{R}^n} \mathrm{e}^{\sqrt{1+|x|^2}}|\phi(x)|\,dx\right)^{1/2}\le C_{\mathrm{init}}\}$. For any $\varepsilon>0$ and $N\in\mathbb{N}$, there exists an FNO $\mathcal{G}_\theta$ such that
\begin{equation}
\begin{aligned}
&\|\G^\dagger-\G_{\theta}\|_{L^2(B_R)}\le\\
&C(3\varepsilon + 3C_{\T{init}}\mathrm{e}^{-\frac{9}{16}R} + (L_{\G}B + C_{\varepsilon, r}^{CC_{\T{ext}}^{(C)}B})L^rN^{-r})
\end{aligned}
\end{equation}
where $C$ depends on $n,r$ and the extension constants, and $L$ is the side length of the torus. The depth and width of the FNO are $\mathcal{O}(\mathrm{depth}(\mathcal{U}^\varepsilon))$ and $\mathcal{O}(N^n\,\mathrm{width}(\mathcal{U}^\varepsilon))$, respectively.
\end{lem}

In the following we use Lemma \ref{sec4.3:lemma1} to show how NOF approximates the optimal CDF; all norms below are $L^2(B_R)$-norms. The time-discretized robust DMZ solution $u_i(t,x)$ approximated by NOF is denoted $\hat u_i(t,x)$ on $[\tau_{i-1},\tau_i]\times B_R$.
\begin{equation}\label{NO_filter_iterate}
\left\{
\begin{aligned}
\hat{u}_i(\tau_i, \cdot) =& \phi_i^{-1}\cdot \G_\theta[\phi_i(x)\cdot \hat{u}_i(\tau_{i-1}, \cdot)]\\
\hat{u}_i(\tau_{i-1}, \cdot) =& \hat{u}_{i-1}(\tau_{i-1}, \cdot),\quad\text{in}\quad[\tau_{i-1}, \tau_i]
\end{aligned}
\right.
\end{equation}
where $\phi_i(x):=\exp(y(\tau_{i-1})\cdot h(x))$. The cumulative error can be computed iteratively. Denoting $e_i:=\hat{u}_i(\tau_i,x)-u_i(\tau_i,x)$, $C_i^\G:=\phi_i^{-1}\circ\G^\dagger\circ\phi_i$, and $e^\G_i:=\phi_i^{-1}\circ(\G_\theta-\G^\dagger)\circ\phi_i\cdot \hat{u}_i(\tau_{i-1},x)$, a direct calculation gives
\begin{equation}
\begin{aligned}
&\|\hat{u}_i - u_i\|_{L^2(B_R)}(\tau_i, \cdot)
\le \|\phi_i^{-1}\circ(\G_\theta-\G^\dagger)\circ\phi_i\cdot \hat{u}_i(\tau_{i-1},\cdot)\|_{L^2(B_R)}\\
&\quad + \|\phi_i^{-1}\circ\G^\dagger\circ\phi_i\big(\hat{u}_i(\tau_{i-1},\cdot)-u_i(\tau_{i-1},\cdot)\big)\|_{L^2(B_R)},
\end{aligned}
\end{equation}
which yields $\|e_i\|_{L^2(B_R)} \le \|e^\G_i\|_{L^2(B_R)} + \|C_i^\G(e_{i-1})\|_{L^2(B_R)}$.

Expanding the second term via \eqref{NO_filter_iterate} and applying the triangle inequality gives
\[
\|C_i^\G(e_{i-1})\|_{L^2(B_R)}
\le \|C_i^\G(e_{i-1}^\G)\|_{L^2(B_R)} + \|C_i^\G\circ C_{i-1}^\G(e_{i-2})\|_{L^2(B_R)}.
\]
Substituting this into the preceding inequality, we get
\begin{equation}
\begin{aligned}
\|e_i\|_{L^2(B_R)} \le& \|e^\G_i \|_{L^2(B_R)} + \|C_i^\G(e_{i-1}^\G)\|_{L^2(B_R)} \\
&+ \|C_i^\G\circ C_{i-1}^\G(e_{i-2})\|_{L^2(B_R)}.
\end{aligned}
\end{equation}
Continuing this process and applying the induction method, we shall have
\begin{equation}\label{ei}
\begin{aligned}
\|e_i\|_{L^2(B_R)} \le& \sum_{k=2}^i\bigg\|C_i^\G\circ C_{i-1}^\G\circ\cdots\circ C_{k+1}^\G(e_k^\G) \bigg\|_{L^2(B_R)} \\
&+ \bigg\|C_i^\G\circ C_{i-1}^\G\circ\cdots\circ C_2^\G(e_1)\bigg\|_{L^2(B_R)}.
\end{aligned}
\end{equation}
We first estimate the norm $\|C_i^\G\circ\cdots\circ C_{k+1}^\G\|_{L^2(B_R)}$ for $i>k$; the proof is given in Appendix \ref{pf:sec4.3_lemma2}.

\begin{lem}\label{sec4.3_lemma2}
Let $i>k$ and $F(x), J(x)$ be defined in Eq.(\ref{dfn_F_J}). 
For a fixed discretization level $N$ and a given truncation level $M>0$, define the event
\[
\begin{aligned}
&\Omega_M := \left\{ \omega : V_0^T(y(\omega)) \le M \right\}, \text{and } \\
 & V_0^T(y):=\sum_{i=0}^N |y(\tau_i)| + \sum_{i=0}^{N-1} |y(\tau_{i+1}) - y(\tau_i)|.
\end{aligned}
\]
Then, for any sample path $y \in \Omega_M$, the following pathwise estimate holds:
\begin{equation}
\begin{aligned}
&\bigg\|C_i^\G\circ \cdots \circ C_{k+1}^\G\bigg\|_{L^2(B_R)}
\le \exp\bigg[M \cdot \|h\|_{\infty} \\
&+ {(i-k)\Delta\tau}\sup_{x\in B_R}\left(-\frac{1}{2}\nabla\cdot F(x) + J(x) + \frac{1}{2}\right)\bigg].
\end{aligned}
\end{equation}
\end{lem}

\begin{rem}
The localization argument does not rely on the specific law of $y$. For a fixed discretization level $N$, the discrete total variation $V_0^T(y)$ defined in Lemma \ref{sec4.3_lemma2} is a finite sum of almost surely finite random variables, so $V_0^T(y)<\infty$ a.s. and $\mathbb{P}(\Omega_M)\to1$ as $M\to\infty$. Thus the error bounds hold with probability arbitrarily close to one for a wide class of observation paths.
\end{rem}

Based on Lemmas \ref{sec4.3_lemma2} and \ref{sec4.3:lemma1}, we can now bound the terms in (\ref{ei}) (proof in Appendix \ref{pf:sec4.3_thm1}).

\begin{thm}\label{sec4.3_thm1}
Let $\{\hat{u}_k\}_{k\ge 0}$ be a constructed density of a NOF on a fixed grid with $N$ time steps. Assume the 
\begin{enumerate}[(1)]
    \item Denote $\tilde{F}(t,x) := -f(x) + \nabla (y(t)\cdot h(x))$ and $\tilde{J}(t, x):= -\nabla\cdot f(x) - \frac12|h(x)|^2 + \frac12\Delta(y(t)\cdot h(x)) - f(x)\cdot \nabla (y(t)\cdot h(x)) + \frac12|\nabla (y(t)\cdot h(x))|^2$ both $C^\infty$ in $x$ and H\"older continuous in $t$ with exponent $\alpha$.
    \item solutions are uniformly bounded: $|u_i(\tau_i,x)|, |u(t,x)| \le C_\rho$ for all $i$ and $x\in B_R$.
\end{enumerate}
For any truncation level $M>0$, define the event
\[
\Omega_M := \left\{ \omega : V_0^T(y(\omega)) \le M \right\}.
\]
Then, on the event $\Omega_M$, the following error estimate holds for each $i \ge 2$: 
\begin{equation}\label{eq:main_error}
\begin{aligned}
&\|\hat{u}_i(\tau_i,\cdot) - u(\tau_i,\cdot)\|_{L^2(B_R)} \\
\le& 
C_{\gamma,i}\, \mathrm{e}^{2M\|h\|_\infty}
\Big( 3\varepsilon + 3C_{\T{init}}\mathrm{e}^{-\frac{9}{16}R} \Big.\\
&\Big.+ (L_{\G}B + C_{\varepsilon, r}^{CC_{\T{ext}}^{(C)}B})L^rN^{-r} \Big) + C_{\alpha,T}\,i^{-\alpha/2},
\end{aligned}
\end{equation}
where
\[
\begin{aligned}
&C_{\gamma,i} := \sum_{k=2}^i \gamma^{i-k} = 
\begin{cases}
\dfrac{1-\gamma^{i-1}}{1-\gamma}, & \gamma \neq 1,\\[1.2ex]
i-1, & \gamma = 1.
\end{cases}, \qquad
\gamma := \mathrm{e}^{\Delta\tau\, \Lambda},\\
&\Lambda := \sup_{x\in B_R}\left(-\frac{1}{2}\nabla\cdot F(x) + J(x) + \frac{1}{2}\right),
\end{aligned}
\]
with $F(x)$ and $J(x)$ as defined in Eq.~(\ref{dfn_F_J}). 
The constants $C_{\text{init}}, L_\G, B, C_{\varepsilon, r}^{CC_{\T{ext}}^{(C)}B}$ are those quantities appearing in Lemma \ref{sec4.3:lemma1}. $C_{\alpha,T}$ is a constant depending on $\alpha, T$.
\end{thm}

\begin{rem}
The final iterative error of NOF has two sources: the approximation error of the forward operator under the numerical YYF scheme, and the discretization error of the robust DMZ solution. Theorem~\ref{sec4.3_thm1} provides a qualitative consistency estimate for a fixed number of observation steps $N$ (hence fixed step size $\Delta\tau=T/N$), with deliberately non-sharp constants. The limit $N\to\infty$ is not of practical interest here, because as $\Delta\tau\to0$ the FKE forward operator $\mathcal{G}^\dagger$ becomes near-identity, a regime where classical single-step solvers are already adequate and operator learning offers no benefit. Thus, only the discretization term $C_{\alpha,T}\,i^{-\alpha/2}$ vanishes with refinement, as expected for a moderate, fixed $N$.
\end{rem}

\section{Numerical experiments}
To demonstrate the reliability and efficiency of the proposed algorithm, we consider two benchmark filtering problems: a one-dimensional cubic filter with nonlinear drift and cubic measurement, and a two-dimensional cubic system with linear drift and cubic observations. Both exhibit severe nonlinearity that challenges conventional filtering methods.

To assess estimation accuracy, we compute the mean-square error (MSE) at each time step and its average (MMSE) over the horizon $[0,T]$: $\text{MSE}(t_k) := \frac{1}{N_{tr}} \sum_{i=1}^{N_{tr}} |\hat{X}^{(i)}(t_k) - X^{(i)}(t_k)|^2$ and $\text{MMSE} := \frac{1}{N} \sum_{k=1}^{N} \text{MSE}(t_k)$, where $\hat{X}^{(i)}(t_k)$ and $X^{(i)}(t_k)$ denote the estimate and the true value at time $t_k$ in the $i$-th trial, and $N_{tr}$ is the number of trials. We also evaluate computational time using the mean time (MT), $\text{MT} = \frac{1}{N_{tr}} \sum_{i=1}^{N_{tr}} T_i$, where $T_i$ denotes the computation time for the $i$-th trial. Throughout all experiments, the number of independent trials is fixed at $N_{tr} = 20$.

\begin{exmp}[1D cubic filter]\label{ex1}
\begin{equation}
\begin{cases}
dX_t = \tanh(X_t)dt + dW_t, \quad\mathbb{E}[(dW_t)^2] = dt,\\
dZ_t = H_0X_t^3dt + dV_t, \quad\mathbb{E}[(dV_t)^2] = sdt,\\
X_0 \sim \sigma_0:=\mathcal{N}(\mu_0 = 0.1,\sigma_0 = 0.05).
\end{cases}
\end{equation}
\end{exmp}
The corresponding forward Kolmogorov equation is expressed as:
\[
\pd{\tilde{u}}{t} = \frac{1}{2}\pdd{\tilde{u}}{x} - \tanh(x)\pd{\tilde{u}}{x} - (\text{sech}^2(x) + \frac{1}{2}\frac{H_0^2x^6}{s})\tilde{u}(t,x)
\]

\begin{table}[!hbt]\small
\centering
\caption{\textbf{Performance comparison of the four methods over 20 independent trials in Example \ref{ex1}.} The bold numbers represent the optimal value among the valid methods; EKF diverges for the cubic observation model, so its MMSE is marked ``--''.}\label{ex1_table}
\begin{tabular}{lcc}
\toprule
Method & MMSE & Time (s) \\
\midrule
NOF(ours) & $0.7132 \pm 0.1100$ & $\mathbf{0.2004 \pm 0.0380}$ \\
PF  & ${0.7034 \pm 0.1150}$ & $1.1851 \pm 0.1000$ \\
LGS  & $0.7195 \pm 0.1130$ & $0.4169 \pm 0.0360$ \\
EKF      & $--$ & $0.002$\\
\bottomrule
\end{tabular}
\end{table}

\begin{figure}[!hbt]
\centering
\includegraphics[width=0.45\linewidth]{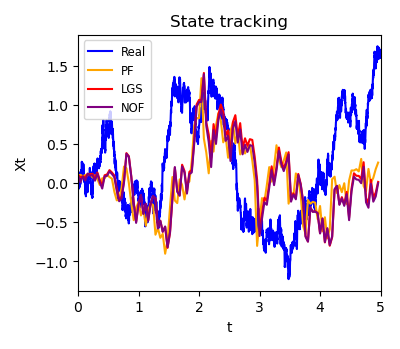}
\includegraphics[width=0.45\linewidth]{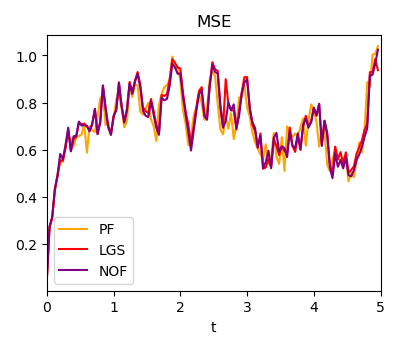}
\caption{\textbf{Performance comparison of state estimation in Example \ref{ex1}}. State track of one Monte-Carlo simulation is shown in the left, where the different color meaning: real state: blue; Particle filter: yellow; Legendre-Galerkin approach: red; Neural operator filter: purple. MMSE obtained by taking 20 repeated Monte-Carlo simulations is presented on the right with the same meaning for colors.}\label{ex1-stateMSE}
\end{figure}

\begin{figure}[!hbt]
    \centering
    \includegraphics[width=0.45\linewidth]{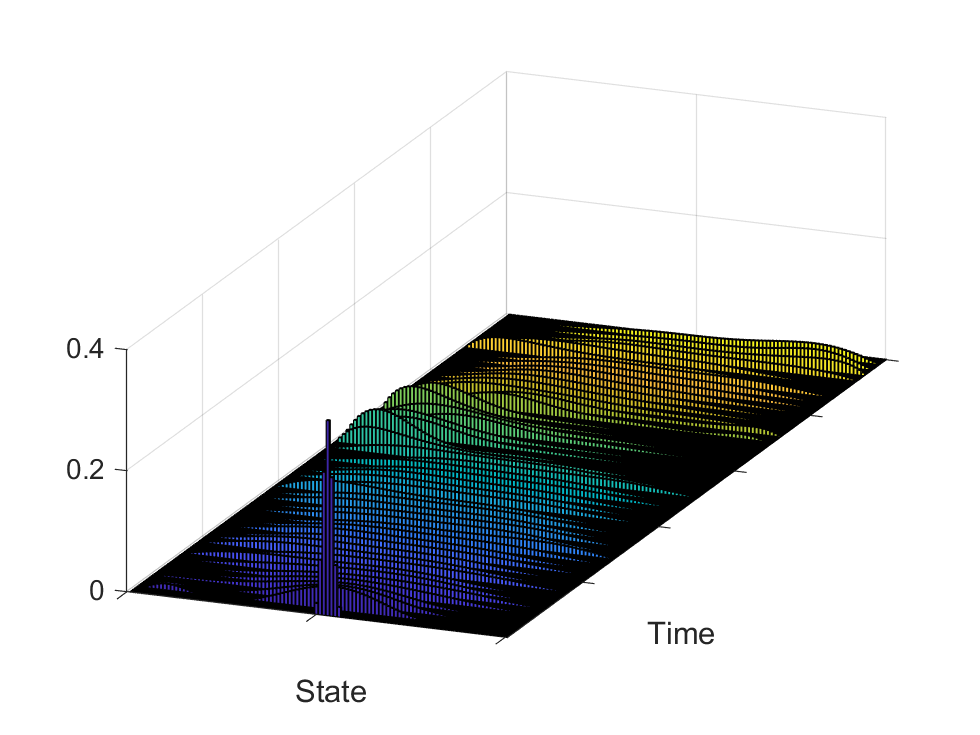}
    \includegraphics[width=0.45\linewidth]{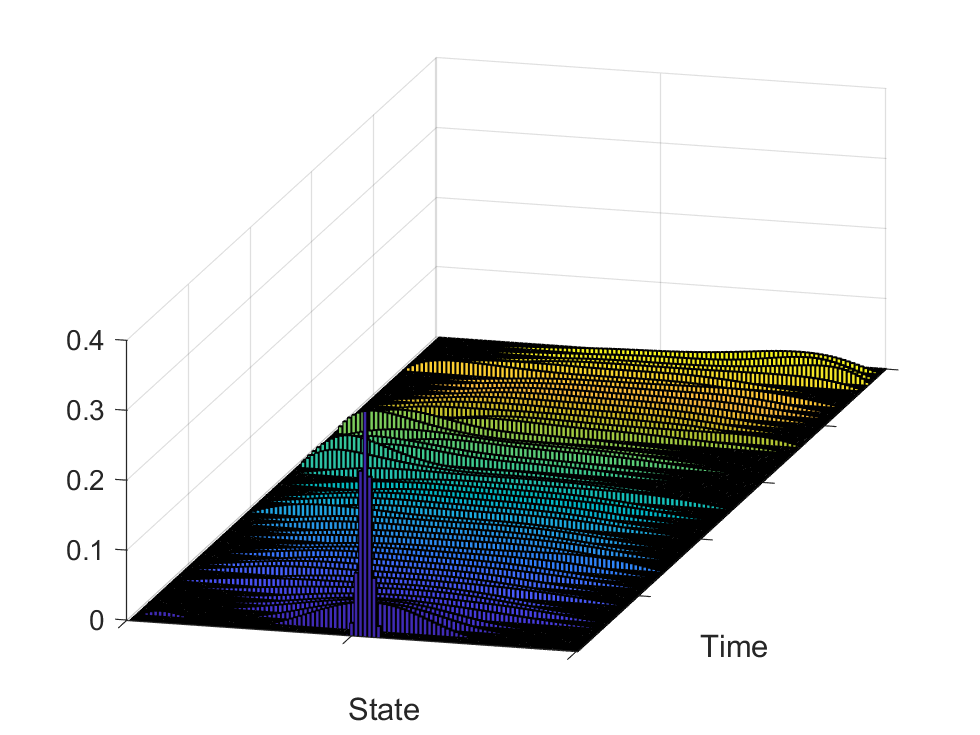}
    \caption{Conditional density function (CDF) computed by Legendre-Galerkin method (left) and neural operator filter (right).}
    \label{fig:ex1_CDF}
\end{figure}
The simulation time is $T=5$ s with observation interval $\Delta t=0.04$; the observation coefficient and noise covariance are $H_0=1.0$ and $s=1.0$, and the initial CDF is $\mathcal{N}(0.1,0.05)$. For the neural operator filter, we sample a training dataset of 235 trajectories on a spatial grid of 100 points, generated by an accurate finite-difference method. We train an operator network with \textit{modes}$=16$, \textit{width}$=40$, and \textit{depth}$=3$ with $\tanh$ activation, using the Adam optimizer with an initial learning rate of $10^{-3}$ that decays by a factor of 0.8 every 100 epochs, for a total of 1000 epochs and a batch size of 300. The particle filter uses 300 particles, and the LGS method uses 8 basis functions. All training is implemented on a single NVIDIA GeForce RTX 2060 Laptop GPU.

Table~\ref{ex1_table} reports the average MMSE and runtime over 20 independent trials. The three valid methods achieve statistically comparable accuracy (MMSE $0.7132\pm0.1100$ for NOF, $0.7034\pm0.1150$ for PF, $0.7195\pm0.1130$ for LGS), with differences well within one standard deviation; EKF diverges for the cubic observation model and is excluded from the accuracy comparison (its MMSE is marked ``--''). In terms of computational efficiency, NOF is significantly faster, consuming only $0.2004\pm0.0380$ seconds per trial, which reduces the runtime by approximately \textbf{83.1\%} compared to PF ($1.1851\pm0.1000$~s) and \textbf{51.9\%} compared to LGS ($0.4169\pm0.0360$~s). The state tracking plot (left panel of Fig.~\ref{ex1-stateMSE}) shows that NOF closely follows the true trajectory, and the time-varying MSE (right panel) confirms that its estimation accuracy remains stable throughout the whole time horizon. These results demonstrate that the neural operator filter not only achieves competitive accuracy but also provides a substantial speedup, making it suitable for real-time applications.

\begin{exmp}[2D cubic system]\label{ex2}
\begin{equation}
\begin{cases}
dX_t = \begin{pmatrix}
a_{11} &a_{12}\\
a_{21} &a_{22}
\end{pmatrix}
\begin{pmatrix}
X_1\\
X_2
\end{pmatrix}
dt + dW_t, \quad\mathbb{E}[dW_tdW_t{^\top}] = I_2dt,\\
dZ_t = \begin{pmatrix}
X_1^3\\
X_2^3
\end{pmatrix}dt + dV_t, \quad\mathbb{E}[dV_tdV_t{^\top}  ] = sIdt,\\ 
X_0 \sim \mathcal{N}([0, 0]^\top, 0.05I_2)
\end{cases}
\end{equation}
\end{exmp}
The corresponding FKE of 2D cubic filter is written as follows.
\begin{equation}\label{}
\begin{aligned}
&\pd{\tilde{u}}{t}(t,x_1,x_2) = \frac{1}{2}(\tilde{u}_{x_1x_1} + \tilde{u}_{x_2x_2}) - (a_{11}x_1+a_{12}x_2)\tilde{u}_{x_1} \\
&- (a_{21}x_1+a_{22}x_2)\tilde{u}_{x_2} - \bigg[(a_{11}+a_{22}) + \frac12(x_1^6+x_2^6)s^{-1}\bigg]\tilde{u},
\end{aligned}
\end{equation}

\begin{figure}[!hbt]
    \centering
    \includegraphics[width=0.8\linewidth]{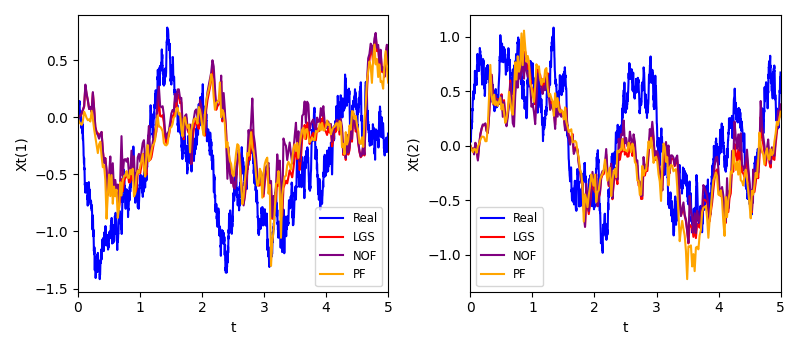}\\
\includegraphics[width=0.8\linewidth]{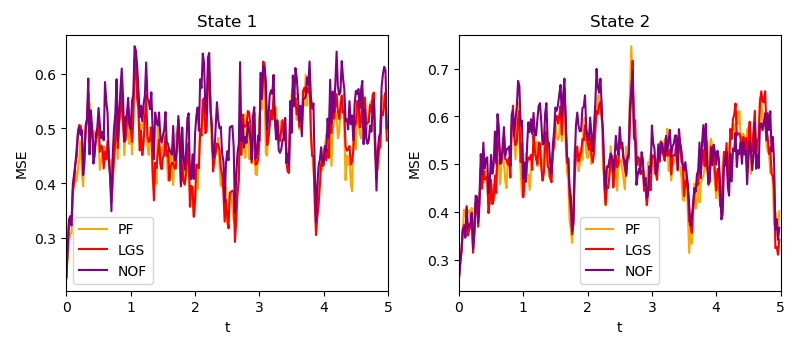}
    \caption{Signal tracking and MSE for 2D cubic system.}
    \label{ex2-state}
\end{figure}

\begin{table}[!hbt]\small
\centering
\caption{\textbf{Performance comparison of four methods over 20 independent trials in Example \ref{ex2}.} The bold numbers represent the optimal value among the valid methods; EKF diverges for the cubic observation model, so its MMSE is marked ``--''.}\label{ex2_table}
\begin{tabular}{lcc}
\toprule
Method & MMSE & Time (s) \\
\midrule
NOF (ours) & $0.5144 \pm 0.0671$ & $\mathbf{1.4116 \pm 0.1065}$ \\
PF         & ${0.4869 \pm 0.0508}$ & $6.4003 \pm 0.2299$ \\
LGS        & $0.4916 \pm 0.0508$ & $3.8479 \pm 0.0793$ \\
EKF & $--$ & $0.01$\\
\bottomrule
\end{tabular}
\end{table}
The observation covariance is set to $S = 0.1 I_2$. The drift coefficients are specified as $a_{11} = -0.3$, $a_{12} = a_{21} = 0.1$, and $a_{22} = 0.2$. For the neural operator filter, we generate a training dataset of 1397 input-output pairs on a spatial resolution of $40\times 40$, and an additional 599 trajectories for testing. The network architecture uses \textit{modes} $=[8,8]$, \textit{width}=40, \textit{depth}=3 and $\tanh$ activation. the relative $L^2$ testing error achieves $0.14\%$. For the particle filter, we use 300 particles to approximate the empirical density in this 2D setting. For the LGS method, the number of basis polynomials is set to 36.

As reported in Table~\ref{ex2_table}, the proposed NOF achieves the lowest computational cost, with an average runtime of $1.4116\pm0.1065$~s, which is approximately \textbf{77.9\%} faster than PF ($6.4003\pm0.2299$~s) and \textbf{63.3\%} faster than LGS ($3.8479\pm0.0793$~s). In terms of estimation accuracy, PF yields the lowest MMSE of $0.4869\pm0.0508$, closely followed by LGS ($0.4916\pm0.0508$) and NOF ($0.5144\pm0.0671$). The state tracking performance for both state dimensions and the corresponding time-varying MSE are depicted in Fig.~\ref{ex2-state}. The figure confirms that NOF closely follows the true trajectory, while its accuracy is competitive with the other methods, and the substantial reduction in runtime makes it a practical choice for real-time filtering tasks.

\section{Conclusion}
This paper proposed a Neural Operator Filter (NOF) that accelerates Yau--Yau nonlinear optimal filtering by removing the heavy on-line burden of solving the Forward Kolmogorov equation (FKE). Using a Fourier Neural Operator (FNO) backbone, NOF learns the prior-to-posterior CDF map directly in an infinite-dimensional function space, avoiding the finite-dimensional truncation of basis-expansion filters and enabling fast, grid-wise distribution updates.

Empirical evaluations on challenging benchmark problems, including the cubic sensor, show that NOF achieves accuracy competitive with strong baselines in mean-square error while substantially improving computational throughput, bringing posterior updates within strict real-time latency constraints. The theoretical analysis establishes stability guarantees and approximation error bounds for both the FKE and the NOF architecture, providing rigorous foundations for the observed performance. Overall, NOF offers a promising operator-learning paradigm for real-time, high-fidelity nonlinear optimal filtering, particularly for complex, multimodal, or heavy-tailed conditional distributions where traditional parametric truncations face inherent limitations.

\section*{Acknowledgments}
This work was supported by grants of the National Natural Science Foundation
of China (No.~12501613, No.~12671549, No.~42450242) and the Tsinghua
University Education Foundation fund.

\bibliographystyle{elsarticle-num}
\bibliography{references}

\appendix
\renewcommand{\thesubsection}{\Alph{section}.\arabic{subsection}}
\section{Proofs}

\subsection{Proof of Proposition \ref{Thm_1}}\label{Proof_Prop4.1}

We derive an \(L^2(B_R)\)-energy estimate for the error function \(w\). 
Taking \(w\) itself as the test function in Eq.~(\ref{Eq::w_R}), multiplying both sides by \(w\), and integrating over \(B_R\), we obtain
\[
\begin{aligned}
	\frac{1}{2}\frac{d}{dt}\|w\|_{L^2(B_R)}^2
	=& \frac{1}{2}\int_{B_R} w\Delta w + \int_{B_R} w F(x)\cdot\nabla w\\
	&+ \int_{B_R} J(x)w^2 + \int_{B_R} w R_i[\tilde v].
\end{aligned}
\]
Using integration by parts together with the homogeneous Dirichlet condition $w=0$ on $\partial B_R$, i.e.,
\[
\int_{B_R}w\Delta w = -\int_{B_R}|\nabla w|^2, \quad
\int_{B_R}wF\cdot\nabla w = -\frac{1}{2}\int_{B_R}w^2\nabla\cdot F
\]
the energy identity becomes
\[
\begin{aligned}
	\frac{1}{2}\frac{d}{dt}\|w\|_{L^2(B_R)}^2
	&= - \frac{1}{2}\int_{B_R}|\nabla w|^2
	- \frac{1}{2}\int_{B_R}w^2\nabla\cdot F \\
	&+ \int_{B_R}J(x)w^2 + \int_{B_R}wR_i[\tilde v].
\end{aligned}
\]
Applying Young's inequality to the last term gives $\int_{B_R}wR_i[\tilde v]\le \frac{1}{2}\int_{B_R}w^2+\frac{1}{2}\int_{B_R}(R_i[\tilde v])^2$.
Consequently,
\[
\begin{aligned}
	&\frac{1}{2}\frac{d}{dt}\|w\|_{L^2(B_R)}^2 + \frac{1}{2}\|\nabla w\|_{L^2(B_R)}^2\\
	&\le \sup_{x\in B_R}\Big[-\frac{1}{2}\nabla\cdot F(x) + J(x) + \frac{1}{2}\Big]\|w\|_{L^2(B_R)}^2\\
	&+ \frac{1}{2}\|R_i[\tilde v]\|_{L^2(B_R)}^2.
\end{aligned}
\]
Multiplying both sides by 2 and then integrating in time interval $[0, T]$, we obtain
\[
\begin{aligned}
	&\|w(T)\|_{L^2(B_R)}^2 + \int_0^T\|\nabla w(s)\|_{L^2(B_R)}^2\,ds\\
	&\le \|w(0)\|_{L^2(B_R)}^2 + \int_0^T\|R_i[\tilde v](s)\|_{L^2(B_R)}^2\,ds
	\\
	&+ \sup_{x\in B_R}\Big[-\nabla\cdot F(x) + 2J(x) + 1\Big]\int_0^T\|w(s)\|_{L^2(B_R)}^2\,ds.
\end{aligned}
\]
Recalling that \(w(0) = \delta\phi\) from the definition of \(w\), this yields precisely the desired estimate \eqref{sec4.1::prop1}.
\hfill\(\square\)

\subsection{Proof of Theorem \ref{Thm_2}}\label{Proof_thm4.1}
 We begin with Eq.(\ref{Diff_form_FKE})
\[
\begin{aligned}
\frac{1}{2}\frac{d}{dt}\int_{B_R}|\tilde u|^2 +& \int_{B_R}\frac{1}{2}|\nabla \tilde u|^2 \\
\le& \sup\bigg[-\frac{1}{2}\nabla\cdot F(x) + J(x) + \frac{1}{2}\bigg]\int_{B_R}\tilde u^2 
\end{aligned}
\]
Let $c:=\sup[-\frac{1}{2}\nabla\cdot F + J(x) + \frac{1}{2}]$. Then
\[
\begin{aligned}
\frac{d}{dt}\left[e^{-2ct}\int_{B_R}\tilde u^2\right] =& 2e^{-2ct}(\frac{1}{2}\frac{d}{dt}\int_{B_R}|\tilde u|^2 - c\int_{B_R}\tilde u^2)\\
\le& -e^{-2ct}\int_{B_R}|\nabla \tilde u|^2 \\
\le& 0
\end{aligned}
\]
Integrating in time over $[0,T]$ gives $\int_{B_R}\tilde u(T,x)^2\,dx\le e^{2cT}\int_{B_R}\tilde u^2(0,x)\,dx$; hence $\G^\dagger$ is Lipschitz with constant $L_\mathcal{G}=\exp\big(T\sup_{x\in B_R}[-\frac{1}{2}\nabla\cdot F+J(x)+\frac{1}{2}]\big)$.

\hfill$\square$
\subsection{Proof of Theorem \ref{thm:FKE_verification}}\label{pf:thm:FKE_verification}

\noindent\textbf{Step 1: Feynman-Kac representation and discretization.}

On the whole space $\mathbb{R}^n$, the Feynman--Kac formula gives
\begin{equation}
u_g(t,x) = \mathbb{E}\left[ \phi(X_t^x) \exp\left( \int_0^t J(X_s^x) \, ds \right) \right],
\label{Pf:A3:EqA2}
\end{equation}
where $X_s^x$ is the solution of $dX_s = F(X_s)ds + dW_s$, $X_0 = x$.

We introduce three successive approximations:

\textit{(i) Initial condition approximation.} For any small $\varepsilon>0$, there is a neural network $\hat{\phi}^\varepsilon$ with $\|\phi-\hat{\phi}^\varepsilon\|_{C^0(B_{R_1})}\le\varepsilon$ for some $R_1>R$ and $\|\hat{\phi}^\varepsilon\|_{C^0(\mathbb{R}^n)}\le C_\phi$; replacing $\phi$ by $\hat{\phi}^\varepsilon$ in \eqref{Pf:A3:EqA2} gives $u_g^\varepsilon$. To control the error we use the \emph{exit event}
\[
\mathcal{E}_{\mathrm{exit}} := \left\{ \max_{0 \le s \le t} \|X_s^x\| \ge R_1 \right\},
\]
where $R_1>R$ is large. On the complement of the exit event the trajectory stays in $B_{R_1}$, where $\hat{\phi}^\varepsilon$ approximates $\phi$ to accuracy $\varepsilon$; on the event itself we use the boundedness of $\phi,\hat{\phi}^\varepsilon$. By Markov's inequality and the standard SDE moment estimate (see \cite{KloedenPlaten1992}), $\mathbb{P}(\mathcal{E}_{\mathrm{exit}})\le C_p(1+R^p)/R_1^p$, so choosing $R_1$ large gives the pathwise bound $|u_g(t,x)-u_g^\varepsilon(t,x)|\le C\varepsilon$ for all $x\in B_R$.

\textit{(ii) Time discretization of the exponential integral.} Discretizing $\int_0^tJ(X_s^x)\,ds$ by the left rectangle rule with $K$ subintervals gives $u_g^{\varepsilon,K}$; splitting on the same exit event and using the local Lipschitz continuity of $J$ on $B_{R_1}$ yields
\begin{equation}
|u_g^{\varepsilon,K}(t,x) - u_g^\varepsilon(t,x)|
\le C\left(\frac{1+R}{\sqrt{K}} + \varepsilon\right), \quad \forall x \in B_R.
\label{Pf:A3:EqA4}
\end{equation}

\noindent\textbf{Step 2: Neural SDE approximation of the solution operator.}

By the universal approximation theorem for Neural SDEs \cite{Kwossek2025}, the solution map $x\mapsto X_t(x)$ can be approximated by a continuous-time Neural SDE; applying the Euler--Maruyama scheme \cite{KloedenPlaten1992} yields a finite-depth ResNet $\mathcal{N}_{k,m}^{(\omega)}(x)$ of the form $\mathcal{N}_{s_{k+1},m}^{(\omega)}(x) = \mathcal{N}_{s_k,m}^{(\omega)}(x) +
\mathcal{F}_\eta(\mathcal{N}_{s_k,m}^{(\omega)}(x)) \Delta s + \Delta W_{k,m}(\omega)$ with initial condition $\mathcal{N}_{0,m}^{(\omega)}(x) = x$ and $\mathcal{F}_\eta$ approximates $F$ and $\Delta W_{k,m}(\omega)$ are Brownian increments; combining the Neural SDE and Euler discretization errors gives $\mathbb{E}\|\mathcal{N}_{K,m}^{(\omega)}(x)-X_t(x)\|\le C\varepsilon$, $\forall x\in B_R$. Now define a Monte Carlo network approximating $u_g^{\varepsilon,K}$ by taking the definition $\mathcal{U}^\varepsilon(\phi, t; \omega)(x) := \frac{1}{M} \sum_{m=1}^M
\hat{\phi}^\varepsilon(\mathcal{N}_{K,m}^{(\omega)}(x))
\cdot \mathcal{E}_{\exp}\left( \sum_{k=0}^{K-1} \mathcal{J}_\delta(\mathcal{N}_{s_k,m}^{(\omega)}(x)) \Delta s \right)$, where $\mathcal{J}_\delta$ approximates $J$ to accuracy $\delta$ and $\mathcal{E}_{\exp}$ approximates $\exp$ to accuracy $\varepsilon_{\exp}$. By the Lipschitz continuity of $\hat{\phi}^\varepsilon,J,\exp$, the Neural SDE error bound above, and the Monte Carlo error $O(1/\sqrt{M})$, we obtain
\begin{equation}
	\begin{aligned}
	&\mathbb{E}\left[|u_g^{\varepsilon,K}(t,x) - \mathcal{U}^\varepsilon(\phi, t; \omega)(x)|\right]\\
	&\le C\left(\varepsilon + \delta + \varepsilon_{\exp} + \varepsilon_\times + \frac{1}{\sqrt{M}}\right),
	\end{aligned}
\label{Pf:A3:EqA7}
\end{equation}
Choosing $M=\lceil\varepsilon^{-2}\rceil$ and $\delta,\varepsilon_{\exp},\varepsilon_\times=O(\varepsilon)$ makes the right-hand side of \eqref{Pf:A3:EqA7} of order $C\varepsilon$.

\noindent\textbf{Step 3: Combining errors and applying the probabilistic method.}

Combining the bounds above with $K=\lceil(1+R)^2/\varepsilon^2\rceil$ gives the pointwise bound
\begin{equation}
\mathbb{E}\left[ |u_g(t,x) - \mathcal{U}^\varepsilon(\phi, t; \omega)(x)| \right]
\le C\varepsilon, \quad \forall x \in B_R.
\label{Pf:A3:EqA9}
\end{equation}

Taking the \(L^q(B_R)\) norm and using Minkowski's inequality (exchanging the expectation and the spatial integral) gives
\begin{equation}
\mathbb{E}\left[ \| u_g(t,\cdot) - \mathcal{U}^\varepsilon(\phi, t; \omega) \|_{L^q(B_R)} \right]
\le C |B_R|^{1/q} \varepsilon.
\label{Pf:A3:EqA10}
\end{equation}
For \(q=2\) this reduces to the \(L^2\) bound; absorbing \(C|B_R|^{1/q}\) into \(\varepsilon\) gives \(\mathbb{E}\|u_g(t,\cdot)-\mathcal{U}^\varepsilon(\phi,t;\omega)\|_{L^q(B_R)}\le\varepsilon\). By the probabilistic method (\cite{de2022generic}, Lemma C.1), some realization $\omega^*$ yields a deterministic network $\mathcal{U}^\varepsilon(\phi,t):=\mathcal{U}^\varepsilon(\phi,t;\omega^*)$ satisfying
\begin{equation}
\| u_g(t,\cdot) - \mathcal{U}^\varepsilon(\phi, t) \|_{L^q(B_R)} \le \varepsilon.
\label{Pf:A3:EqA12}
\end{equation}
For \(q=2\) this is the \(L^2(B_R)\) bound. Since $\mathcal{U}^\varepsilon$ is a composition of smooth $\tanh$ networks, its derivatives up to order $\ell$ are bounded on the compact domain $B_R$.
\hfill$\square$

\subsection{Proof of Corollary \ref{Cor:FKE_verification}}\label{pf:Cor:FKE_verification}

By employing the truncation error estimate from \cite{yau2008real}, for sufficiently large $R>0$ and $0 \le t \le T$,
\begin{equation}
\int_{B_R} |u_g(t,x) - u(t,x)| \, dx
\le C_{n,T} e^{-\frac{9}{8}R} \int_{\mathbb{R}^n} e^{\sqrt{1+|x|^2}} |\phi(x)| \, dx.
\end{equation}

Since $u_g$ and $u$ are uniformly bounded on $[0,T]\times B_R$ by the maximum principle (with $J\le\bar J$), $\|u_g(t,\cdot)-u(t,\cdot)\|_{L^\infty(B_R)}\le C_\infty$ for some $C_\infty>0$.
Since $\|u_g-u\|_{L^1(B_R)}\le C_{n,T}e^{-\frac{9}{8}R}\int_{\mathbb{R}^n} e^{\sqrt{1+|x|^2}}|\phi|\,dx$ and $\|u_g-u\|_{L^\infty(B_R)}\le C_\infty$, the interpolation inequality $\|\cdot\|_{L^2}^2\le\|\cdot\|_{L^1}\|\cdot\|_{L^\infty}$ yields
\begin{equation}
\|u_g(t,\cdot)-u(t,\cdot)\|_{L^2(B_R)}\le C_{n,T}\, e^{-\frac{9}{16}R}\left(\int_{\mathbb{R}^n} e^{\sqrt{1+|x|^2}}|\phi(x)|\,dx\right)^{1/2}
\label{Pf:A3:EqA30}
\end{equation}

The triangle inequality applied to \eqref{Pf:A3:EqA12} and \eqref{Pf:A3:EqA30} gives the desired bound with the factor $e^{-\frac{9}{16}R}$.

\hfill$\square$

\subsection{Proof of Lemma \ref{sec4.3:lemma1}}\label{pf:sec4.3:lemma1}
We first recall the approximation result for the pseudo-spectral projection on a torus.

\begin{lem}[\cite{Shen2011Spectral}]\label{lem:scaled_interpolation}
Let $n\in\mathbb{N}$, $r>n/2$, and $\mathbb{T}_L^n:=\mathbb{R}^n/(L\mathbb{Z}^n)$ the $n$-dimensional torus of side length $L>0$. For $N\in\mathbb{N}$, let $Q_N\circ\mathcal{E}_N$ be the trigonometric interpolation (pseudo-spectral projection) operator onto frequencies of modulus at most $N$. Then, for every $f\in H^r(\mathbb{T}_L^n)$ and $0\le k\le r$,
\begin{equation}
\| f - (Q_N \circ \mathcal{E}_N) f \|_{H^k(\mathbb{T}_L^n)}
\le C(r,n) \, L^{\,r-k} \, N^{-(r-k)} \, \| f \|_{H^r(\mathbb{T}_L^n)},
\end{equation}
where the constant $C(r,n) > 0$ is independent of $L$, $N$, and $f$.
\end{lem}

\textit{Proof of the Lemma \ref{sec4.3:lemma1}}.
The total error decomposes into three parts:
\begin{equation}\label{pf:A5:Eq1}
\begin{aligned}
&\|\G^\dagger - \G_{\theta}\|_{L^2(B_R)} \le \|\G^\dagger - \mathcal{U}^\varepsilon(\phi, t)\|_{L^2(B_R)}\\
&+ \|\mathcal{U}^\varepsilon(\phi, t)-\hat{\G}\|_{L^2(B_R)} + \|\hat{\G}-\G_\theta\|_{L^2(B_R)}
\end{aligned}
\end{equation}
By Corollary \ref{Cor:FKE_verification}, for any $\varepsilon>0$ and sufficiently large $R$, the network $\mathcal{U}^\varepsilon(\phi,t)$ approximates $\G^\dagger$ with accuracy $\varepsilon+C_{\mathrm{trunc}}(R)\le\varepsilon+C_{n,T}C_{\T{init}}e^{-\frac{9}{16}R}$, bounding the first term. 

For the second term, we first verify the stability of the network $\mathcal{U}^\varepsilon(\phi,t)$. For any $\phi_1,\phi_2\in C^r(B_R)\cap S_\phi$,
\[
\begin{aligned}
&\|\mathcal{U}^\varepsilon(\phi_1, t) - \mathcal{U}^\varepsilon(\phi_2, t)\|_{L^2(B_R)}\\
&\le\|\mathcal{U}^\varepsilon(\phi_1, t) - \tilde{u}(t, x;\phi_1) \| + \|\tilde{u}(t, x;\phi_1) - \tilde{u}(t, x;\phi_2) \|\\
&+ \|\tilde{u}(t, x;\phi_2) - \mathcal{U}^\varepsilon(\phi_2, t)\|_{L^2(B_R)}\\
&\le2\varepsilon+2C_{\T{init}}C_{n, T}e^{-\frac{9}{16}R}+L_\G\|\phi_1-\phi_2\|_{L^2(B_R)}
\end{aligned}
\]
Taking $\phi_1=\phi$ and $\phi_2=\mathcal{R}\circ P_N\circ\mathcal{E}(\phi)$ in the stability bound gives $\T{(II)}\le 2\varepsilon+2C_{\T{init}}C_{n,T}e^{-\frac{9}{16}R}+L_\G\|\phi-\mathcal{R}\circ P_N\circ\mathcal{E}(\phi)\|_{L^2(B_R)}$.
For the last term,
\[
\begin{aligned}
&\|\phi-\mathcal{R}\circ P_N\circ \mathcal{E}(\phi)\|_{L^2(B_R)}\le \|\mathcal{E}\phi - P_N(\mathcal{E}\phi)\|_{L^2(\mathbb{T}^n_{L})}\\
&\le C(r,n)C_{\T{ext}}^{(H)} L^rN^{-r}\|\phi\|_{H^r(B_R)}\le C(r,n)C_{\T{ext}}^{(H)}B L^rN^{-r}
\end{aligned}
\]
where the second inequality uses Lemma \ref{lem:scaled_interpolation} and the stability of $\mathcal{E}$, and the third uses $C^r(B_R)\hookrightarrow H^r(B_R)$.

The third term is estimated similarly: $\|\hat\G-\G_\theta\|_{L^2(B_R)}\le C(r,n)C_{\T{ext}}^{(H)}L^rN^{-r}\|\hat\G\|_{H^r(B_R)}$, where $\|\hat\G\|_{H^r(B_R)}\le C_{\varepsilon,r}^{CC_{\T{ext}}^{(H)}B}$ because $\|\mathcal{R}\circ P_N\circ\mathcal{E}(\phi)\|_{C^r(B_R)}\le C\,C_{\T{ext}}^{(C)}\|\phi\|_{C^r(B_R)}\\
\le C\,C_{\T{ext}}^{(C)}B$.

Combining the three estimates yields the bound in the statement.
\hfill$\square$
\subsection{Proof of Lemma \ref{sec4.3_lemma2}}\label{pf:sec4.3_lemma2}

Recall the $L^2(B_R)$-operator norm bound for the FKE forward map $\G^\dagger$: by Proposition \ref{Thm_2},
\begin{equation}\label{Lemma4.1_eq1}
\|\G^\dagger\|_{L^2(B_R)} 
\le \exp\!\Big(\Delta\tau \sup_{x\in B_R}\big[-\tfrac12\nabla\cdot F(x) + J(x) + \tfrac12\big]\Big).
\end{equation}
Now consider the product of compositions \(C_i^\G\circ \cdots \circ C_{k+1}^\G\) for \(i>k\). 
Using the submultiplicativity of the operator norm and the explicit form of the \(C_j^\G\)'s, 
we obtain
\[
\begin{aligned}
&\|C_i^\G\circ \cdots \circ C_{k+1}^\G\|_{L^2(B_R)}\\
&\le \|\phi_i^{-1}\|_{L^2(B_R)}
   \Big(\prod_{j=k+2}^{i} \|\phi_j\circ\phi_{j-1}^{-1}\|_{L^2(B_R)}\Big)
   \|\phi_{k+1}\|_{L^2(B_R)}\\
&\times \|\G^\dagger\|_{L^2(B_R)}^{i-k} \\
&= \Big(\|e^{-y(\tau_{i-1})h}\|_{L^2(B_R)}\Big.\\
  &\Big.\times\prod_{j=k+1}^{i-1} \|e^{(y(\tau_{j+1})-y(\tau_j))h}\|_{L^2(B_R)}
   \|e^{y(\tau_k)h}\|_{L^2(B_R)}\Big)
   \|\G^\dagger\|_{L^2(B_R)}^{i-k}
\end{aligned}
\]
Each multiplication operator \(e^{g h}\) (where \(g\) is a scalar function of time) has \(L^2(B_R)\)-norm bounded by \(\exp(\|g\|_\infty \|h\|_\infty)\). Hence
\[
\begin{aligned}
&\|e^{-y(\tau_{i-1})h}\|_{L^2(B_R)}
   \prod_{j=k+1}^{i-1} \|e^{(y(\tau_{j+1})-y(\tau_j))h}\|_{L^2(B_R)}
   \|e^{y(\tau_k)h}\|_{L^2(B_R)} \\
&\le \exp\!\Big(\big(|y(\tau_{i-1})| + \sum_{j=k+1}^{i-1} |y(\tau_{j+1})-y(\tau_j)| + |y(\tau_k)|\big)\|h\|_\infty\Big).
\end{aligned}
\]
Combining this with \eqref{Lemma4.1_eq1} yields the final estimate.

\hfill$\square$

\subsection{Proof of Theorem \ref{sec4.3_thm1}.}\label{pf:sec4.3_thm1}
We start from the error decomposition $\|\hat{u}_i(\tau_i,\cdot) - u(\tau_i,\cdot)\|_{L^2(B_R)} 
\le \|\hat{u}_i(\tau_i,\cdot) - u_i(\tau_i,\cdot)\|_{L^2(B_R)}+ \|u_i(\tau_i,\cdot) - u(\tau_i,\cdot)\|_{L^2(B_R)}$, where $u_i$ is the numerical Yau--Yau solution using the exact $\mathcal{G}^\dagger$; the second term is bounded by $C_{\alpha,T}i^{-\alpha/2}$ via Theorem~C of \cite{yau2008real}.

For the first term, the recursive formula (\ref{ei}) gives $\|\hat{u}_i(\tau_i,\cdot) - u_i(\tau_i,\cdot)\|_{L^2(B_R)}
\le \sum_{k=2}^i \big\|C_i^\mathcal{G} \circ \cdots \circ C_{k+1}^\mathcal{G}(e_k^\mathcal{G})\big\|_{L^2(B_R)}\\
+ \big\|C_i^\mathcal{G} \circ \cdots \circ C_2^\mathcal{G}(e_1)\big\|_{L^2(B_R)}$. By Lemma~\ref{sec4.3_lemma2} and the bound on $\|e_k^\mathcal{G}\|$, the first term has an upper bound $\sum_{k=2}^i \big\|C_i^\mathcal{G} \circ \cdots \circ C_{k+1}^\mathcal{G}(e_k^\mathcal{G})\big\|_{L^2(B_R)}\le e^{2M\|h\|_\infty} \sum_{k=2}^i \gamma^{\,i-k} \|\mathcal{G}_\theta - \mathcal{G}^\dagger\|_{L^2(B_R)} = C_{\gamma,i}\, e^{2M\|h\|_\infty} \|\mathcal{G}_\theta - \mathcal{G}^\dagger\|_{L^2(B_R)}$, where $e^{2M\|h\|_\infty}$ is a crude uniform bound absorbing the multiplication operators $\phi_k^{\pm1}$ in the error injections $e_k^\mathcal{G}$ together with the composition norms of Lemma~\ref{sec4.3_lemma2}; it is not meant to be sharp. For the initial error term, $\big\|C_i^\mathcal{G} \circ \cdots \circ C_2^\mathcal{G}(e_1)\big\|_{L^2(B_R)}
\le e^{M\|h\|_\infty + T\Lambda} \|\mathcal{G}_\theta - \mathcal{G}^\dagger\|_{L^2(B_R)}$. Substituting Lemma \ref{sec4.3:lemma1} and combining all terms yields \eqref{eq:main_error}.
\hfill$\square$

\end{document}